\documentclass[UTF-8,11pt]{article}
\usepackage[T1]{fontenc}
\usepackage{latexsym,amssymb,amsmath,amsfonts,amsthm}
\usepackage{graphics}
\usepackage[demo]{graphicx}
\usepackage[section]{placeins}\usepackage{graphicx}
\usepackage{mathrsfs}
\usepackage{subfigure}
\usepackage{float}
\usepackage{color}

\usepackage{tikz}

\usetikzlibrary{calc}
\usetikzlibrary{intersections,through}
\usepackage{pgfplots}
\usepgfplotslibrary{patchplots}
\usepgfplotslibrary{colormaps}
\usepackage{mathtools}

\usepackage{epstopdf}
\newcommand{\R}{{\mathbb R}}

\newcommand{\be}{\begin{eqnarray}}
\newcommand{\ben}{\begin{eqnarray*}}
\newcommand{\en}{\end{eqnarray}}
\newcommand{\enn}{\end{eqnarray*}}

\newtheorem{theorem}{Theorem}[section]

\newtheorem{lemma}[theorem]{Lemma}

\newtheorem{remark}[theorem]{Remark}

\newtheorem{algorithm}{Algorithm}[section]

\definecolor{rot}{rgb}{0,0,0}
\definecolor{hw}{rgb}{0,0,0}

\definecolor{mgq}{rgb}{0,0,0}
\newcommand{\mgq}{\color{mgq}}
\begin{document}
\renewcommand{\theequation}{\arabic{section}.\arabic{equation}}

\title{Imaging a moving point source in $\R^3$ from the time of arrival at sparse observation points}

	\author{Guanqiu Ma\footnotemark[3], Haonan Zhang\footnotemark[2], and Hongxia Guo\footnotemark[1] }
	
	\date{}
	\maketitle
	
	\renewcommand{\thefootnote}{\fnsymbol{footnote}}
	\footnotetext[3] {School of Mathematical Sciences, Sichuan Normal University, Chengdu, 610068, China (gqma@sicnu.edu.cn).}
	
	\footnotetext[1]{Corresponding author: School of Mathematical Sciences, and Institute of Mathematics and Interdisciplinary Sciences, Tianjin Normal University, Tianjin, 300387, China(hxguo@tjnu.edu.cn).}

    \footnotetext[2] {Beijing Computational Science Research Center, Beijing, 100193, China (zhanghaonan@csrc.ac.cn).}
	\renewcommand{\thefootnote}{\arabic{footnote}}

	\begin{abstract}

In this paper, we introduce a novel numerical method for reconstructing the trajectory within three-dimensional space, where both the emission moment and spatial location of the point source are unknown. Our approach relies solely on measuring the time of arrival at five or seven properly chosen observation points.  By utilizing the distinctive geometric configuration of these five or seven observation points, we  establish the uniqueness of the trajectory and emission moment of the point source through rigorous mathematical proofs. Moreover, we analyze the stability of our proposed method. The effectiveness of the method is also verified by numerical experiments.

\vspace{.2in} {\bf Keywords}: {\bf inverse moving source problem, wave equation, uniqueness, time of arrival}
	\end{abstract}

	\section{Introduction}

In a homogeneous and isotropic medium occupying the whole space $\mathbb{R}^3$ with constant sound speed $c > 0$,  we consider the dynamic wave propagation generated by a moving point source. The source follows a trajectory $a(t) \in C(t_{\min}, t_{\max}) $, with $0 < t_{\min} < t_{\max}$ and 
travels at a subsonic speed, satisfying
\begin{equation}\label{speed-p}
|a'(t)| < c \quad \mbox{for all} \quad  t \in [t_{\min}, t_{\max}].
\end{equation}
During the observation period, the source emits a finite sequence of pulses described by
	\begin{equation}
		S(x,t) = \sum\limits_{j=1}^{J} \delta(x-a(t)) \delta(t-t_j) \ell(t), \quad t_{\min} \leq t_1 < \cdots < t_J \leq t_{\max},
	\end{equation}
	where $\delta$ denotes the Dirac delta function, and $\ell(t)$ is a real-valued continuous function obeying the   coercivity condition $|\ell(t)| \geq \ell_0 >0$.
	
The radiated acoustic field $U(x,t)$ satisfies the initial value problem 
	\begin{equation}
		\left\{
		\begin{aligned}
			&c^{-2}\frac{\partial^2 U}{\partial t^2} = \Delta U + S(x,t), \quad &&(x, t) \in \mathbb{R}^3 \times \mathbb{R}^+, \mathbb{R}^+ \coloneqq \{t\in \R: t>0\},\\
			&U(x,0)=\partial_t U(x,0) = 0, &&x\in \mathbb{R}^3.
		\end{aligned}
		\right.
	\end{equation}

By convolving the fundamental solution $G(x;t,c) = \frac{\delta(t-c^{-1}|x|)}{4\pi |x|}$ with the source term, the explicit solution reads
	\begin{align*}
		U(x,t)& = G(x;t,c) * S(x,t) \coloneqq \int_{\mathbb{R}^+} \int_{\mathbb{R}^3}
		G(x-y; t-\tau,c)S(y,\tau)\,dyd\tau\\
&=\sum\limits_{j=1}^{J} \frac{\delta (t-t_j -c^{-1}|x-a(t_j)|)}{4\pi |x-a(t_j)|} \ell(t_j), \quad x \notin \Gamma,
	\end{align*}
where $\Gamma \coloneqq \{x: x=a(t),\, t\in [t_{\min}, t_{\max}]\}$.
	For a given receiver location $x_k$, the signal is observed only at discrete {\bf Times of Arrival}:
	\begin{equation}\label{toa}
	T_{jk} = t_j + c^{-1}|x_k-a(t_j)|,\quad j=1,\cdots,J.
	\end{equation}

	In this paper we are interested in the following inverse problem:
	\begin{description}
	\item[(IP):] Recovery the trajectory $\Gamma$ and emission moments $t_j,\, j=1,\cdots,J$ using the Time of Arrival
	$$\{T_{jk}: j=1,2,\cdots,J; \, k=1,\cdots,K\}$$
	of the observation points $\{x_k: k=1,\cdots,K\}$, $x_k \notin \Gamma$.
\end{description}
	
The problem of reconstructing the location and emission moments of a moving point source arises widely in radar positioning, navigation, communication transmission, sonar and other practical applications.
{\mgq Numerical methods for inverting moving point sources based on  arrival time data can generally be categorized into two appoaches: the Time of Arrival (TOA) method and the Time Difference of Arrival (TDOA) method. There requires not only the recording of signal arrival times at observation points but also prior knowledge of the signal's emission time in the TOA method \cite{FKK1999, S1999}. Theoretically, by obtaining emission and arrival time data from four non-coplanar observation points, a system of equations can be constructed to invert the position of the source at the emission moment. However, this method imposes stringent requirements on time synchronization between the source and observation points, as calibration errors can degrade inversion accuracy. In contrast, the TDOA method eliminates the need for time synchronization between the source and observation points, offering distinct advantages in ensuring precision.
Existing research on TDOA \cite{CC2003, CNAS2014, GM2015, GG2003, LQS2012, WSZG2019, DZY2020, TOY2021, TWC2022} demonstrates that scholars have understood the theoretical approach of determining a hyperboloid using arrival time differences from two observation points. Additionally, algorithms for constructing equation systems based on arrival time differences and solving for source locations via least squares methods have been extensively studied. Nevertheless, current studies exhibit deficiencies in the observation point configurations and the rigorous proof of solution uniqueness in inversion. To address these gaps, this paper is concerned with the optimal measurement geometry design and the rigorous mathematical proof of the uniqueness of the inversion solution.

In this paper, we address these limitations by proposing  two novel TDOA-based observation point deployment schemes:
\begin{itemize}
 \item The five-point positioning scheme, which achieves point source inversion in three-dimensional space with a minimal number of observation points but presents certain limitations in practical applications.
 \item The seven-point positioning scheme, which requires only seven observation points to stably invert both the emission time and spatial location of an arbitrary moving point source.
     \end{itemize}
For both configurations, we provide rigorous uniqueness theorems and stability analyses, offering solid theoretical support for practical implementation. }

The remainder of the paper is organized as follows. Section \ref{uniq} presents the uniqueness results for the proposed inverse problem. In section \ref{alg}, we introduce details of the reconstruction algorithm and its stability. Numerical experiments validating our method are reported in Section \ref{num}. Finally, Section \ref{con} concludes the paper and outlines possible directions for future work.

\section{Uniqueness}\label{uniq}

The ordering of the arrival times is determined under the assumption \eqref{speed-p}.
\begin{lemma}\label{time-order}
	Let $x_k$ be a fixed point satisfying  the assumption \eqref{speed-p}, the arrival times satisfy
	\begin{equation*}
	T_{jk} > T_{lk} \quad \mbox{if and only if} \quad j>l.
	\end{equation*}
\end{lemma}

\begin{proof}
	Suppose, for contradiction, that  $j>l$ but $T_{jk} \leq T_{lk}$. By the triangular inequality and the mean value theorem, we have
	\begin{equation*}
		\begin{aligned}
			T_{lk}-T_{jk} &= t_l -t_j +c^{-1}(|x_k-a(t_l)|-|x_k-a(t_j)|) \\
			& \leq t_l -t_j +c^{-1} |a(t_j) - a(t_l)| \\
			& = (t_j -t_l)(c^{-1}|a'(t)|-1)\\
			& < 0,
		\end{aligned}
	\end{equation*}
for some $t\in (t_l,t_j)$.  This contradicts our hypothesis $T_{jk} \leq T_{lk}$.
\end{proof}

Owing to the ordering of the arrival times and the independence of each emission, the uniqueness of the point source at each emission moment can be established.

There are five observation points denoted by ${x_1, x_2, \cdots, x_5 }$. From \eqref{toa}, it follows that 
\begin{equation}
	\left(a(t_j)-x_k\right)^2 = c^2(T_{jk}^2 -2 T_{jk} t_j +t_j^2),\quad j=1,\cdots, J, \,k=1,\cdots,5.
\end{equation}
Assume without loss of generality that $x_1=0$(the origin). {\mgq  Subtracting the  equations corresponding to $x_1$ from those for the remaining points yields }
\begin{equation}
	x_k \cdot a(t_j) + c^2(T_{j1}-T_{jk}) t_j= \frac{c^2}{2} (T_{j1}^2 - T_{jk}^2) + \frac{x_k^2}{2}, \quad j=1,\cdots, J, \,k=2,\cdots,5.
\end{equation}
This system can be rewritten in matrix form as
\begin{equation}
	A_{j} X_{j} = b_{j}, \quad j=1,\cdots, J,
\end{equation}
where
\begin{equation*}
	X_{j} := (a(t_j) \quad t_j)^{\top},\quad b_{j} (k-1) := \frac{c^2}{2} (T_{j1}^2 - T_{jk}^2) + \frac{x_k^2}{2},\, k=2,\cdots,5
\end{equation*}
and
\begin{equation*}
	A_{j}(k-1,:) := (x_k \quad c^2(T_{j1}- T_{jk})),\, k=2,\cdots,5.
\end{equation*}

\begin{theorem}\label{un-5}(Uniqueness with five observation points)
	If $A_{j}$ is invertible, then the locations of the point source $a(t_j)$ and the emission moments $t_j$ can be determined from the  the time-of-arrival data $T_{jk}$, $k=1,2,\cdots,5$.
\end{theorem}

\begin{proof}
	If $A_{j}$ is invertible, then $ X_{j} = A_{j}^{-1} b_{j}$ uniquely determines $a(t_j)$ and $t_j$.
\end{proof}

In practice, however the invertiblity of matrix $A_{j}$ in Theorem \ref{un-5} is not guaranteed. To address this, we consider an alternative configuration of observation points that is more practical.
{\mgq Since a hyperboloid can be uniquely determined by three colinear observation points, we distribute seven observation points along the three coordinate axes to ensure the uniqueness of the point source.}
We next demonstrate the uniqueness of both the source location $a(t_j)$  and the emission time $t_j$ under this configuration.

There are seven observation points denoted by ${x_1, x_2, \cdots, x_7 }$. Specifically, $x_1, x_2, x_3$ are collinear on line $l_1$, $x_1, x_4, x_5$ are collinear on line $l_2$, and $x_1, x_6, x_7$ are collinear on line $l_3$. Moreover, $x_6$ and $x_7$ are symmetric with respect to $x_1$, that is, $|x_1 - x_6| = |x_1 - x_7|$. In addition, the lines $l_1$, $l_2$, and $l_3$ are mutually perpendicular in pairs (see Figure \ref{seven-points}).
	
	\begin{figure}[!ht]
		\centering
	\begin{minipage}{0.4\textwidth}
		\centering
		
		\scalebox{0.6}{\begin{tikzpicture}
		\draw[->] (-4,0) -- (4,0) node[above] {$y(l_1)$} coordinate(x axis);
		\draw[->] (3,3) -- (-3,-3) node[right] {$x(l_2)$} coordinate(y axis);
		\draw[->] (0,-4) -- (0,4) node[right] {$z(l_3)$} coordinate(z axis);
		\draw (0.3,0) node [below] {$x_1$};
		\draw (2,0) node [below] {$x_2$};
		\draw (-3,0) node [below] {$x_3$};
		\draw (2.6,2.5) node [below] {$x_4$};
		\draw (1.6,1.5) node [below] {$x_5$};
		\draw (0,3) node [left] {$x_6$};
		\draw (0,-3) node [left] {$x_7$};		
		\fill (0,0) circle (2pt);
		\fill (2,0) circle (2pt);
		\fill (-3,0) circle (2pt);
		\fill (2.5,2.5) circle (2pt);	
		\fill (1.5,1.5) circle (2pt);	
		\fill (0,3) circle (2pt);	
		\fill (0,-3) circle (2pt);							
		\end{tikzpicture} }
		\caption{Seven points in $\mathbb{R}^3$.}
		\label{seven-points}	
	\end{minipage}
	\begin{minipage}{0.4\textwidth}
		\centering
 		
		\scalebox{0.6}{\begin{tikzpicture}
			\draw[->] (-4,0) -- (4,0) node[above] {$y$} coordinate(x axis);
			\draw[->] (0,-4) -- (0,4) node[right] {$z$} coordinate(z axis);
			\draw [thick,blue,smooth] (3,2) .. controls (1,0) .. (3,-2);
			\draw [thick,,red,smooth] (4,2) .. controls (-0.5,0) .. (4,-2);
			\draw (0.3,0) node [below] {$x_1$};
			\draw (2,0) node [below] {$x_2$};
			\draw (-3,0) node [below] {$x_3$};
			\fill (0,0) circle (1pt);
			\fill (2,0) circle (1pt);
			\fill (-3,0) circle (1pt);
			\fill (2.2,1.2) circle (2pt);
			\fill (2.2,-1.2) circle (2pt);
			\draw [dotted] (2.2,1.2) -- (2.2,-1.2);
			\coordinate (A) at (2.2,0);
			\draw (A) rectangle ($(A)+(.1,.1)$);
		\end{tikzpicture} }
		\caption{The intersection of two hyperboloids.(The red line shows the equation \eqref{s32} and the blue line shows the equation \eqref{s12}.)}
		\label{x1x2x3}
		
	\end{minipage}
\end{figure}

\begin{theorem}(Uniqueness with seven observation points)
	The locations of the point source $a(t_j)$ and the emission moments $t_j$ can be determined by the time of arrival $T_{jk}$, $k=1,2,\cdots,7$.
\end{theorem}

\begin{proof}

From the arrival times $T_{j1},T_{j2}$ and $T_{j3}$, we obtain two independent equations for the source location. Without loss of generality, assume  $T_{j2} < T_{j1} < T_{j3}$ for a fixed $j$. Then,
\begin{equation}\label{s32}
	|a(t_j) - x_3| - |a(t_j) - x_2| = c(T_{j3} - T_{j2}),
\end{equation}
and
\begin{equation}\label{s12}
	|a(t_j) - x_1| - |a(t_j) - x_2| = c(T_{j1} - T_{j2}).
\end{equation}
 Given that $x_1$, $x_2$, and $x_3$ are colinear along  the $y$-axis, each equation defines one branch of a two-sheet hyperboloid of revolution about the $y$-axis.  Their cross-sectional diagrams in the $yOz$ plane can be observed in Figure \ref{x1x2x3}. The intersection of two hyperboloids can result in, at most, a circle perpendicular to the $y$-axis.
Similarly, the arrival times $T_{j1}$, $T_{j4}$, and $T_{j5}$ yield another hyperbolic system with respect to the $x$-axis, giving rise to a second circle orthogonal to the $x$-axis.
The intersection of these two circles determines, at most, two candidate source locations due to the symmetry of the $xOy$ plane. The additional arrival times $T_{j6}$ and $T_{j7}$ are used to resolve this ambiguity. 
For example, if $T_{j6} < T_{j7}$, then $a(t_j)$ must lie closer to $x_6$, thus selecting the unique physical solution. 

Notably, when the source lies on a coordinate axis, four distinct arrival times corresponding to three collinear observation points and a single point off-axis are sufficient to uniquely determine $a(t_j)$. Finally, once $a(t_j)$ is known, the emission time $t_j$ is uniquely given by 
\[ t_j = T_{j1} - c^{-1} |a(t_j) - x_1|. \]

This completes the proof. 
\end{proof}

	\section{Algorithm and stability}\label{alg}

In this section, we present the algorithm to calculate the locations of the point source $a(t_j)$ and the emission moments $t_j$, $j=1,\cdots,J$.

\begin{algorithm} \label{algo5}(Five observation points)
	\\
	{\bf Step 1:} From $T_{11}$, $T_{12}$, $T_{13}$, $T_{14}$ and $T_{15}$, we build $A_{1}$ and $b_{1}$ by
	\begin{equation*}
		A_{1}(k-1,:) := (x_k \quad c^2(T_{11}- T_{1k})),\,b_{1} (k-1) := \frac{c^2}{2} (T_{11}^2 - T_{1k}^2) + \frac{x_k^2}{2},\, k=2,\cdots,5.
	\end{equation*}
	\\
	{\bf Step 2:} Calculate $X_{1}$ from
	\begin{equation}\label{X-com}
		X_{1} =(a(t_1) \quad t_1)^{\top} = (A_{1})^{-1} b_{1}.
	\end{equation}
	\\
	{\bf Step 3:} Repeat the above two steps and calculate $\{a(t_j);t_j\}$ from $\{T_{j1},T_{j2},T_{j3},T_{j4},T_{j5}\}$, $j=2,\cdots,J$.
\end{algorithm}

\begin{algorithm} \label{algo7}(Improved TDOA with seven observation points)
	\\
	{\bf Step 1:} From $T_{11}$, $T_{12}$, $T_{13}$, $T_{14}$,$T_{15}$,$T_{16}$ and $T_{17}$, we build $A^{(1)}$ and $b^{(1)}$ by
	\begin{equation*}
		A_{1}(k-1,:) := (x_k \quad c^2(T_{11}- T_{1k})),\,b_{1} (k-1) := \frac{c^2}{2} (T_{11}^2 - T_{1k}^2) + \frac{x_k^2}{2},\, k=2,\cdots,7.
	\end{equation*}
	\\
	{\bf Step 2:} Calculate $X_{1}$ from
	\begin{equation}\label{X-com}
		X_{1} =(a(t_1) \quad t_1)^{\top} = (A_{1}^{\top}A_{1})^{-1} A_{1}^{\top} b_{1}.
	\end{equation}
	\\
	{\bf Step 3:} Repeat the above two steps and calculate $\{a(t_j);t_j\}$ from $\{T_{j1},T_{j2},T_{j3},T_{j4},T_{j5},T_{j6},T_{j7}\}$, $j=2,\cdots,J$.
\end{algorithm}

Next, we will analyze the stability of $a(t_j)$ and $t_j$ in the above algorithms. In what follows, we shall use the 2-norm.

\begin{theorem}\label{sta5}
	Let $A_{j}$ is invertible. Assume that $\tilde{T}_{jk}$ is the measured time of arrival with perturbations, and that $\tilde{a}(t_j)$ and $\tilde{t}_j$ are the emission moments and locations of the point source calculated from the perturbed data. There exists a uniform constant $ M>0$, such as $||\delta A|| := ||\tilde{A}_{j} -A_{j}||\leq Mc^2\delta$ and $||\delta b|| := ||\tilde{b}_{j} -b_{j}||\leq Mc^2\epsilon$ for $j=1,\cdots,J$. Then, the following estimate holds
	\begin{equation}
		||\tilde{X}_{j} - X_{j}|| \leq {\mgq \frac{Mc^2(\epsilon +||A^{-1}||\,||b||\delta)}{||A^{-1}||^{-1}-Mc^2\delta} },\quad j=1,2,\cdots,J.
	\end{equation}
\end{theorem}

\begin{proof}
	Without loss of generality, we use the case $j=1$ as an example to prove, and denote that $A := A_{1}, \delta A := \tilde{A}_{1} -A_{1}, b := b_{1}, \delta b := \tilde{b}_{1} -b_{1}, X := X_{1}$ and $\delta X := \tilde{X}_{1} - X_{1}$.

	Thus, we have
	\begin{equation*}
		(A+\delta A)(X+ \delta X) = b+\delta b,
	\end{equation*}
	which means that
	\begin{equation*}
		\delta X = A^{-1}(\delta b - \delta A X -\delta A  \delta X).
	\end{equation*}
	Using the triangle inequality, we have $$||\delta X|| \leq ||A^{-1}|| \, (||\delta b|| + ||\delta A|| \, ||X|| +||\delta A|| \, ||\delta X||).$$
	From {\mgq $||X||\leq ||A^{-1}||||b||$}, $||\delta A|| \leq Mc^2\delta$ and $||\delta b|| \leq Mc^2\epsilon$, we get
	{\mgq \begin{equation*}
		||\delta X|| \leq \frac{||A^{-1}||\left(||\delta b||+||\delta A|| \, ||x||\right)}{1- ||A^{-1}||\,||\delta A||} \leq \frac{Mc^2(\epsilon +||A^{-1}||\,||b||\delta)}{||A^{-1}||^{-1}-Mc^2\delta}.
	\end{equation*} }
\end{proof}

\begin{theorem}\label{sta5.2}
	If $A_{j}$ is singular, we assume that $\tilde{T}_{jk}$ is the measured time of arrival with perturbations, and that $\tilde{a}(t_j)$ and $\tilde{t}_j$ are the emission moments and locations of the point source calculated from the perturbed data. There exists a uniform constant $ M>0$, such as $||\delta A|| := ||\tilde{A}_{j} -A_{j}||\leq Mc^2\delta$ and $||\delta b|| := ||\tilde{b}_{j} -b_{j}||\leq Mc^2\epsilon$ for $j=1,\cdots,J$. Then, the following estimate holds
	\begin{equation}
		||\tilde{X}_{j} - X_{j}|| \leq {\mgq\frac{Mc^2 (||A||+Mc^2\delta) (\epsilon + \delta ||(A^{\top}A)^{-1}||\,||A||\,||b||)}{||(A^{\top}A)^{-1}||^{-1}-2||A||Mc^2\delta-(Mc^2\delta)^2}, }\quad j=1,2,\cdots,J.
	\end{equation}
\end{theorem}

\begin{proof}
	Without loss of generality, we use the case $j=1$ as an example to prove, and denote that $A := A_{1}, \delta A := \tilde{A}_{1} -A_{1}, b := b_{1}, \delta b := \tilde{b}_{1} -b_{1}, X := X_{1}$ and $\delta X := \tilde{X}_{1} - X_{1}$.

	Thus, we have
	\begin{equation*}
		(A+\delta A)^{\top}(A+\delta A)(X+ \delta X) = (A+\delta A)^{\top} (b+\delta b),
	\end{equation*}
	which means that
	{\mgq 
	\begin{equation*}
		\begin{aligned}
			\delta X &= (A^{\top}A)^{-1} [A^{\top}\delta b + (\delta A)^{\top} \delta b- A^{\top}\delta A X -(\delta A)^{\top}\delta AX - \\
			& A^{\top}\delta A  \delta X - (\delta A)^{\top}A \delta X - (\delta A)^{\top}\delta A \delta X].
		\end{aligned}
	\end{equation*}}
	Using the triangle inequality, we have 
	{\mgq $$||\delta X|| \leq ||(A^{\top}A)^{-1}|| \, \left[(||A||+||\delta A||)(||\delta b|| + ||\delta A||\,||X||)+(2||A||+||\delta A||)||\delta A||\,||\delta X|| \right].$$}
	From {\mgq $||X||\leq ||(A^{\top}A)^{-1}||\,||A|| \,||b||$}, $||\delta A|| \leq Mc^2\delta$ and $||\delta b|| \leq Mc^2\epsilon$, we get
	{\mgq \begin{equation*}
			||\delta X|| \leq \frac{(||A||+||\delta A||)(||\delta b|| + ||\delta A||\,||X||)}{||(A^{\top}A)^{-1}||^{-1} - ||\delta A||^2 -2||\delta A||\,||A||} 
			\leq \frac{Mc^2 (||A||+Mc^2\delta) (\epsilon + \delta ||(A^{\top}A)^{-1}||\,||A||\,||b||)}{||(A^{\top}A)^{-1}||^{-1}-2||A||Mc^2\delta-(Mc^2\delta)^2}.
	\end{equation*} }
\end{proof}

\section{Numerical Examples}\label{num}

In this section, we present two numerical examples to validate the proposed algorithm in three-dimensional space. In both cases, we consider the effect of measurement perturbations and compare the numerical reconstructions against the theoretical stability estimates.

Using Algorithm \ref{algo7}, the trajectory of a moving point source is recovered from the recorded arrival times at seven strategically placed observation points: $x_1=(0,0,0)$, $x_2=(3,0,0)$, $x_3=(-3,0,0)$, $x_4=(0,3,0)$, $x_5=(0,-3,0)$, $x_6=(0,0,3)$, and $x_7=(0,0,-3)$. The wave speed is set to $c=1$. The arrival moments at each sensor are denoted by $T_{jk}$ for $k=1,\dots,7$.

\begin{figure}[htp]
\centering
\subfigure[$a_r$]{
\includegraphics[scale=0.27]{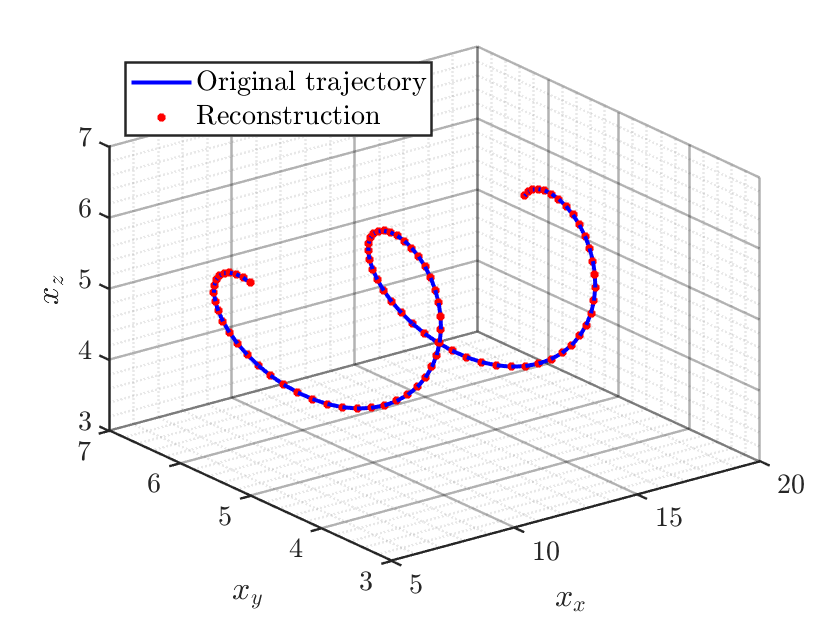}
}
\subfigure[$\mathrm{error}$]{
\includegraphics[scale=0.25]{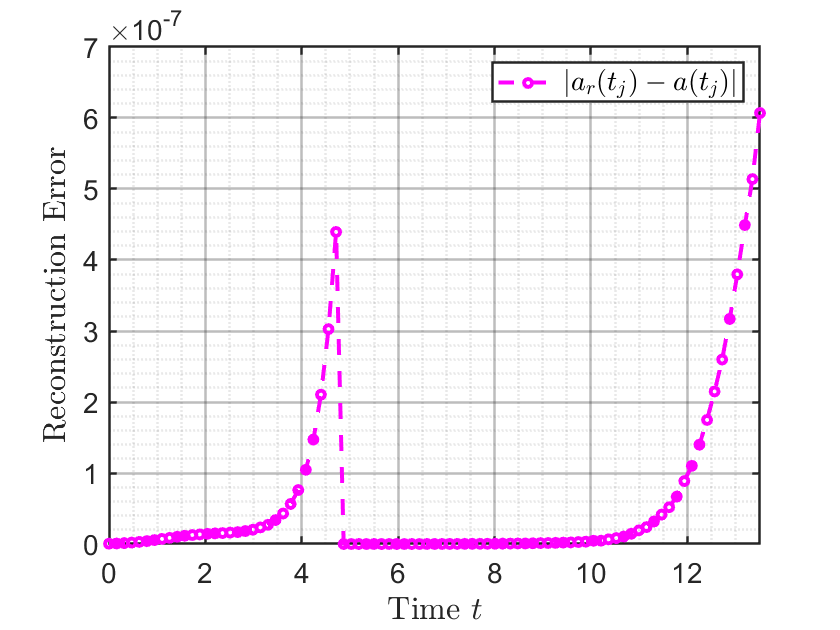}
}

\subfigure[$a_r$ under noise]{
\includegraphics[scale=0.27]{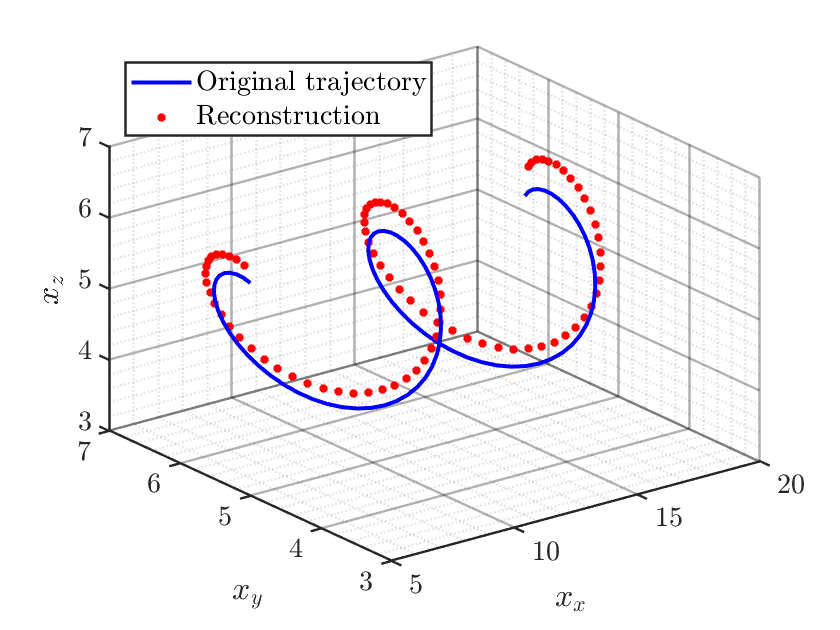}
}
\subfigure[error under noise]{
\includegraphics[scale=0.25]{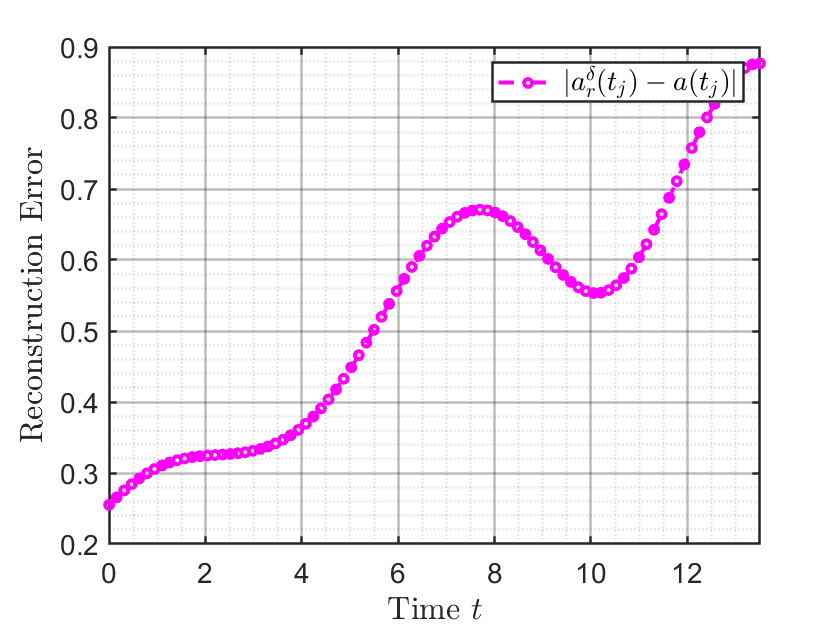}
}
\caption{Reconstruction of a spiral trajectory using seven observation points. (a) Reconstructed trajectory without noise; (b) corresponding reconstruction error; (c) reconstructed trajectory with $1\%$ noise; (d) corresponding reconstruction error under noise.}
\label{fig:1}
\end{figure}

\textbf{Example 1.} Consider a point source moving along a spiral trajectory defined by $a(t) = (t+5, \sin t + 5, \cos t + 5)$ for $t \in [0, 4\pi + 1]$. The time discretization is given by $t_j =  (j-1)\frac{\pi}{20}$ for $j=1,\dots,86$. Figure \ref{fig:1} shows the reconstructed trajectory and corresponding errors.

Figures \ref{fig:1}(a) and \ref{fig:1}(b) illustrate the reconstruction $a_r(t_j)$ and the reconstruction error $\mathrm{error}(t_j) = \|a_r(t_j) - a(t_j)\|$ in the absence of noise ($\delta=0$). The results confirm that the trajectory is accurately recovered, with error below $10^{-6}$. Figures \ref{fig:1}(c) and \ref{fig:1}(d) show the reconstruction $a_r^{\delta}(t_j)$ and the error ${error}^\delta(t_j) = \|a_r^{\delta}(t_j) - a(t_j)\|$ for a noise level of $\delta = 1\%$. While the noise leads to noticeable perturbations, the reconstructed trajectory remains a good approximation of the true path, demonstrating the algorithms robustness.

\begin{figure}[htp]
\centering
\subfigure[$a_r$]{
\includegraphics[scale=0.27]{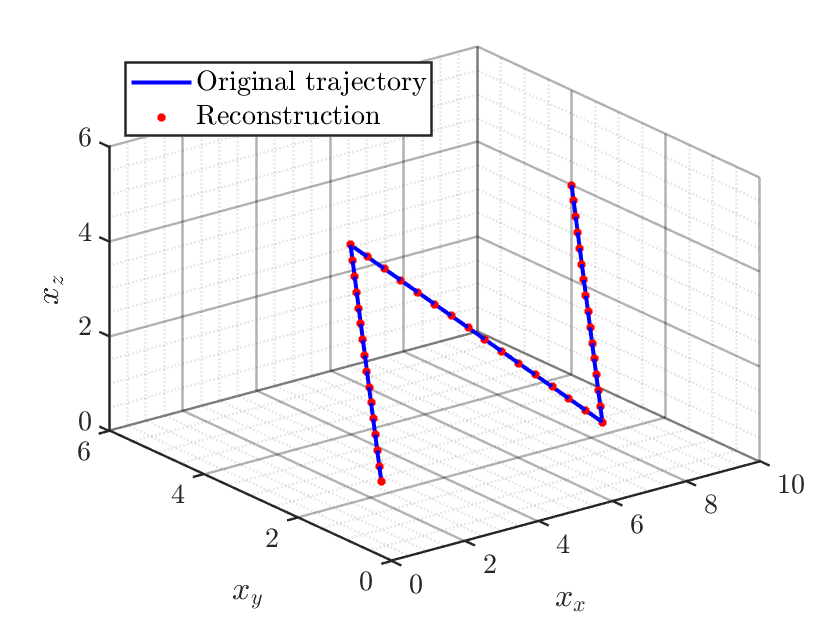}
}
\subfigure[$\mathrm{error}$]{
\includegraphics[scale=0.25]{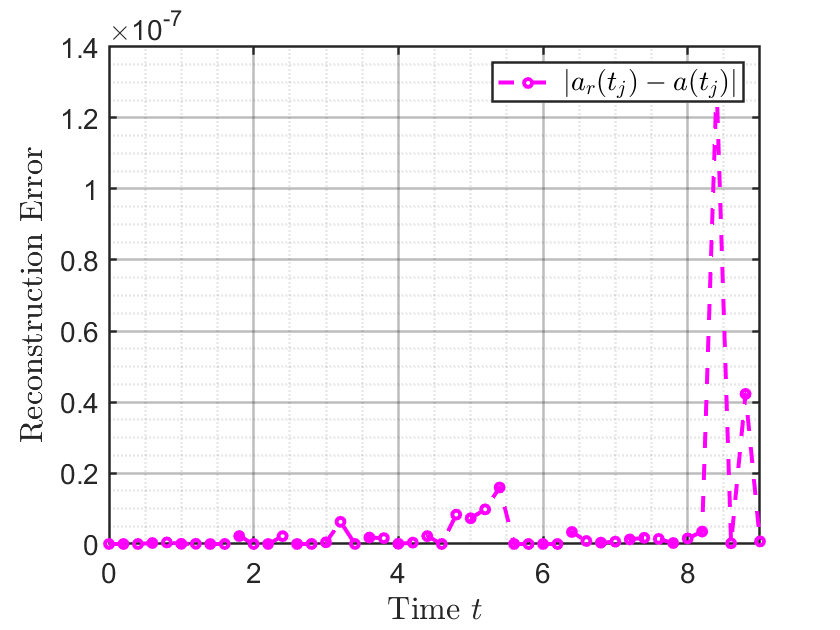}
}

\subfigure[$a_r$ under noise]{
\includegraphics[scale=0.27]{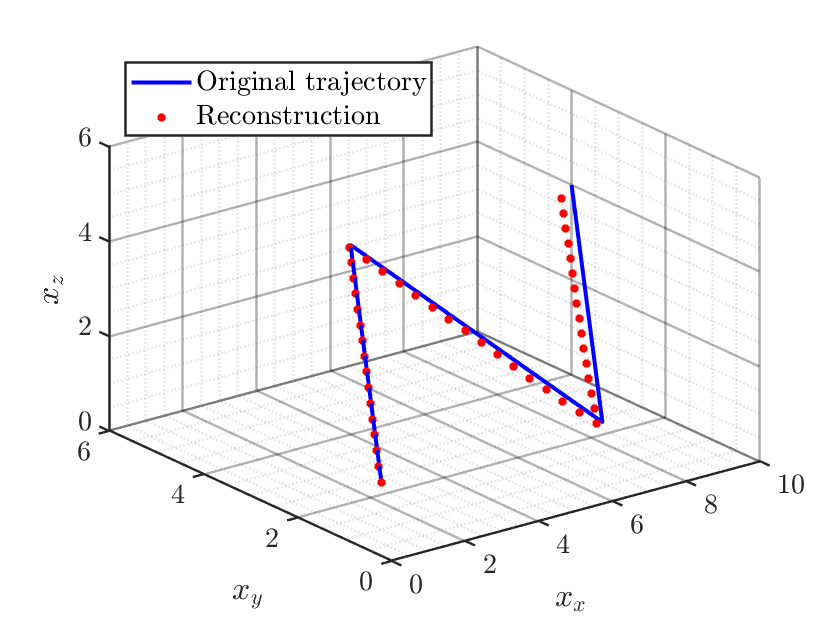}
}
\subfigure[error under noise]{
\includegraphics[scale=0.25]{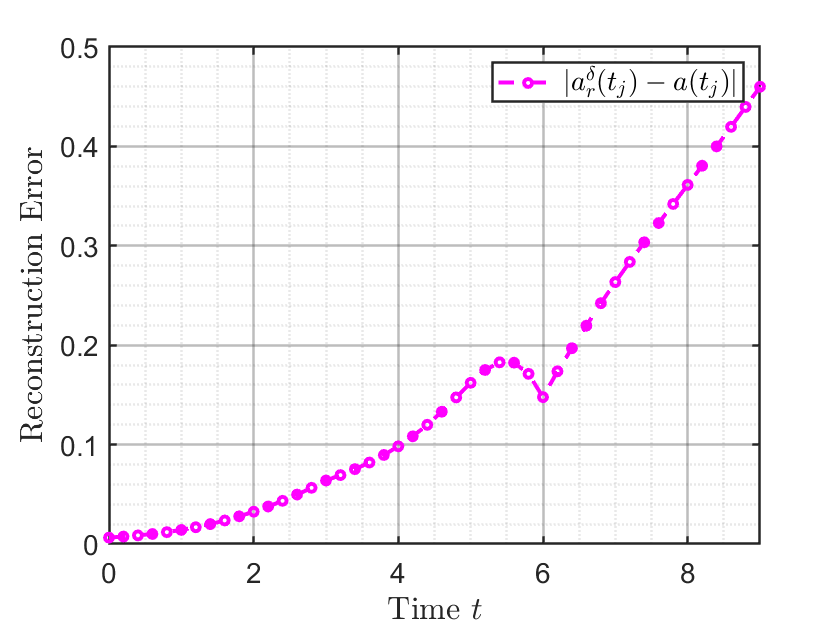}
}
\caption{Reconstruction of a piecewise linear (folded) trajectory using seven observation points. (a) Reconstructed trajectory without noise; (b) corresponding reconstruction error; (c) reconstructed trajectory with $1\%$ noise; (d) corresponding reconstruction error under noise.}

\label{fig:2}
\end{figure}

\textbf{Example 2.} Next, we consider a piecewise linear (folded) trajectory:
\begin{equation*}
a(t) = \begin{cases}
[1 + t,\; 4 + t,\; 7 + t], & 0 \leq t \leq 3, \\
[-2 + t,\; 7 - t,\; -2 + t], & 3 < t \leq 6, \\
[-5 + t,\; 10 + t,\; -5 + t], & 6 < t \leq 9.
\end{cases}
\end{equation*}
The trajectory is discretized at $t_j = (j-1) \times 0.2$ for $j=1,\dots,46$. Figure \ref{fig:2} presents the reconstructions and errors with and without noise.

Figures \ref{fig:2}(a) and \ref{fig:2}(b) demonstrate accurate recovery of the fold line in the noise-free case, with reconstruction errors again below $10^{-6}$. When perturbations are introduced ($\delta = 1\%$), as shown in Figures \ref{fig:2}(c) and \ref{fig:2}(d), the trajectory is slightly disturbed, but the general structure is still well preserved. This confirms the algorithms resilience against moderate noise.

\begin{figure}[htp]\label{fig:stability}
\centering
\includegraphics[scale=0.3]{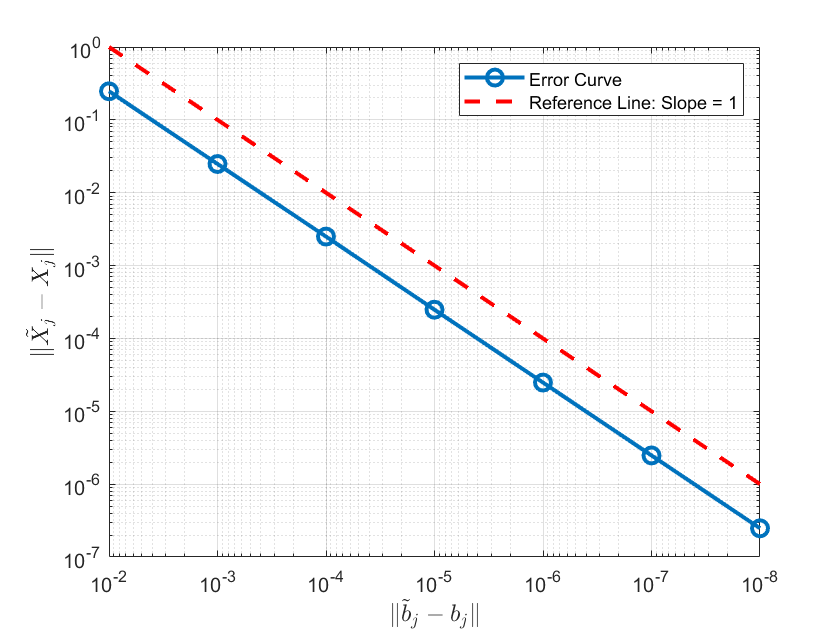}
\caption{Numerical verification of the stability estimate in Theorem \ref{sta5}. The plot shows the reconstruction error $\|\tilde{X}_j - X_j\|$ versus the data perturbation $\|\tilde{b}_j - b_j\|$ on a logarithmic scale. }
\end{figure}

\textbf{Stability Verification.} To validate the theoretical stability estimate stated in Theorem \ref{sta5.2}, we conduct a numerical test shown in Figure \ref{fig:stability}. Theorem \ref{sta5.2} asserts that the reconstruction error satisfies

	\begin{equation*}
		||\tilde{X}_{j} - X_{j}|| \leq  {\mgq\frac{Mc^2 (||A||+Mc^2\delta) (\epsilon + \delta ||(A^{\top}A)^{-1}||\,||A||\,||b||)}{||(A^{\top}A)^{-1}||^{-1}-2||A||Mc^2\delta-(Mc^2\delta)^2}, \quad j=1,2,\cdots,J.}
	\end{equation*}
indicating first-order stability with respect to data perturbations.

Figure \ref{fig:stability} displays the reconstruction error $\|\tilde{X}_j - X_j\|$ versus noise level $\|\tilde{b}_j - b_j\|$ on a log-log scale. The results reveal a linear trend with slope $1$, consistent with the theoretical prediction. A reference line with slope $1$ is also included for comparison. The close match between the numerical results and the theoretical line confirms the first-order convergence and validates the robustness of the algorithm.
\begin{remark}
A separate numerical verification of Theorem \ref{sta5} is not included here, as the stability result shares the same first-order behavior as that of Theorem 3.1, which has already been confirmed numerically. Given the similarity in the structure of the error estimates, we refrain from presenting redundant numerical experiments.
\end{remark}


\section{Conclusions and future work}\label{con}

After the number and position of observation points are specially designed, the improved TDOA algorithm can effectively recover the location and emission moment of the moving point source. The stability analysis also shows that the error of this algorithm depends on the wave velocity $c$. When the wave velocity $c$ is large enough, the error of the location of the point source may be significant, even when the measurement data error is small.
The time of arrival about the point source cannot be measured effectively in $\R^2$, and this improved TDOA algorithm is no longer feasible since there is no Huygens principle in two-dimensional space. If there are multiple moving point sources in the space, it is difficult to separate the time of arrival corresponding to multiple point sources, which is a problem we will consider in future.

\section*{Acknowledgments}

The work of H. Guo is partially supported by Natural Science Foundation of Tianjin Municipal Science and Technology Commission (No.24JCQNJC00850).

\end{document}